\documentclass[12pt]{article}
\usepackage{amsmath}
\usepackage{amssymb}
\usepackage{tabularx}
\usepackage{enumerate}
\usepackage{graphicx}
\usepackage{color}

\newtheorem{thm}{Theorem}[section]
\newtheorem{lem}[thm]{Lemma}
\newtheorem{cor}[thm]{Corollary}

\newtheorem{rmk}{Remark}[section]

\newtheorem{pppp}{Proof}

\newcommand{\qed}{\hspace{1em}\mbox{\raisebox{0.65ex}{\fbox{}}}}

\numberwithin{equation}{section}

\newcommand{\be}{\begin{equation}}
\newcommand{\ee}{\end{equation}}
\newcommand\bes{\begin{eqnarray}} \newcommand\ees{\end{eqnarray}}
\newcommand{\bess}{\begin{eqnarray*}}
\newcommand{\eess}{\end{eqnarray*}}

\newcommand{\R}{\mathbb{R}}
\newcommand{\bpf}{{\bf Proof:\ \ }}
\newcommand{\epf}{\mbox{}\hfill $\Box$}

\begin{document}

\thispagestyle{empty}

\title{Reproduction numbers and the expanding fronts for a diffusion-advection SIS model in heterogeneous time-periodic
environment\thanks{The work is partially supported by the NSFC of China (Grant No. 11371311 and 11501494), the High-End Talent Plan of Yangzhou University.}}
\date{\empty}

\author{Jing Ge$^{a}$, Chengxia Lei$^{b}$, Zhigui Lin$^{a,}$\thanks{Corresponding author. Email address: zglin68@hotmail.com (Z. Lin).} \\
{\small $^a$School of Mathematical Science, Yangzhou University, Yangzhou 225002, China}\\
{\small $^b$Department of Mathematics, University of Science and Technology of China, Hefei 230026, China}
}

 \maketitle

\begin{quote}
\noindent
{\bf Abstract.} { 
This paper deals with a simplified SIS model, which describes the transmission of the disease in time-periodic heterogeneous environment. To understand the impact of spatial heterogeneity of environment and small advection on the persistence and eradication of an infectious disease, the left and right free boundaries are introduced to represent the expanding fronts. The basic reproduction numbers $R_0^D$ and $R_0^F(t)$, which depends on spatial heterogeneity, temporal periodicity and advection, is introduced. A spreading-vanishing dichotomy is established and sufficient conditions for the spreading and vanishing of the disease are given. The asymptotic spreading speeds for the left and right fronts are also presented.
 }


\end{quote}

\section{Introduction}

To understand the transmission of infectious diseases, many mathematical models have been made and investigated.
Considering spatial diffusion and environmental heterogeneity, Allen, Bolker, Lou and Nevai in \cite{AL} proposed
an SIS epidemic reaction-diffusion model
\begin{eqnarray}
\left\{
\begin{array}{lll}
S_{t}-d_S\Delta S=-\frac{\beta (x) SI}{S+I}+\gamma(x)I,\; &  x\in\Omega, \ t>0, \\
I_{t}-d_I\Delta I=\frac{\beta (x) SI}{S+I}-\gamma(x)I,\; &   x\in \Omega,\ t>0,\\
\frac {\partial S}{\partial \eta}=\frac {\partial I}{\partial \eta}=0,\;&  x\in \partial \Omega, \ t>0,
\end{array} \right.
\label{Aa1}
\end{eqnarray}
where $S(x, t)$ and $I(x, t)$ represent the density of susceptible and infected individuals at
location $x$ and time $t$, respectively, the positive constants $d_S$ and $d_I$ denote the corresponding diffusion rates for
the susceptible and infected populations, $\beta(x)$ and $\gamma(x)$ are positive H$\ddot{o}$lder continuous functions,
 which account for spatial dependent rates of disease contact transmission and disease recovery at $x$, respectively.
 The term $\frac{\beta (x) SI}{S+I}$ is the standard incidence of disease.

 As in \cite{AL},
we say that $x$ is a {\bf high-risk site} if the local disease transmission rate $\beta(x)$ is greater
than the local disease recovery rate $\gamma(x)$. An {\bf low-risk site} is defined in a similar manner.
The habitat $\Omega$ is characterized as {\bf high-risk}
( or {\bf low-risk} ) if the spatial average $(\frac{1}{|\Omega|}\int_{\Omega}\beta(x)dx)$
of the transmission rate is greater than ( or less than ) the spatial
average $(\frac{1}{|\Omega|}\int_{\Omega}\gamma(x)dx)$ of the recovery rate, respectively.

In some recent work \cite{Pe1, Pe3, PZ}, Peng et al. further investigated the asymptotic
behavior and global stability of the endemic equilibrium for system (\ref{Aa1}) subject to the Neumann boundary conditions,
and provided much understanding of the impacts of large and small diffusion rates of the susceptible
and infected population on the persistence and extinction of the
disease.

To focus on the new phenomena induced by spatial heterogeneity of environment, we assume that the population $N(x,t)$ is constant in space for all time, that is, $N(x,t)\equiv N^*$ for $x\in \Omega$ and $t\geq 0$ and consider he corresponding free boundary problem
\begin{eqnarray}
\left\{
\begin{array}{lll}
I_{t}-d_I I_{xx}+\alpha I_{x}=(\beta (x,t) -\gamma(x,t))I-\frac {\beta (x,t)}{N^*} I^2,\; &g(t)<x<h(t),\; t>0,   \\
I(g(t),t)=0,\, g'(t)=-\mu I_{x}(g(t),t),&t>0,\\
I(h(t),t)=0,\, h'(t)=-\mu I_{x}(h(t),t),  & t>0, \\
g(0)=-h_0,\, h(0)=h_0,\, I(x,0)=I_{0}(x), & -h_0\leq x\leq h_0,
\end{array} \right.
\label{a3}
\end{eqnarray}
where $x=g(t)$ and $x=h(t)$ are the moving left and right
boundaries to be defined,  $h_0,\, d_I,\, \alpha$ and $\mu $ are positive constants.
$\alpha$ and $\mu$ are referred as the advection rate and the expanding capability, respectively.
$\beta(x,t), \gamma (x,t)\in C^{\nu_0, \frac{ \nu_0} 2}(\R \times [0,\infty))$ for some $\nu_0\in (0,1)$,
 which account for spatial dependent rates of disease contact transmission and disease recovery, respectively. We assume that
$\beta(x,t)$ and $\gamma (x,t)$ are positive and bounded, that is, there exist positive constants $\beta_1, \beta_2, \gamma_1$ and $\gamma_2$ such that $\beta_1\leq \beta(x,t)\leq\beta_2$ and $\gamma_1\leq\gamma (x,t)\leq \gamma_2$ in $\R \times [0,\infty)$.
Considering periodic environment, we assume that $\beta(x,t), \gamma(x,t)$
 are periodic in $t$ with the same period $T$ (i.e., $\beta(x,t+T)=\beta (x,t)$, $\gamma(x,t+T)=\gamma(x, t)$ for all $t\in \R$).
 Further, in the paper we assume
$$(H_1)\;\lim_{x\rightarrow\pm \infty}\beta(x,t)=\beta_{\infty}(t),\;
\lim_{x\rightarrow \pm \infty}\gamma(x,t)=\gamma_{\infty}(t)\, \textrm{uniformly for}\, t\in [0, T],$$
which means that far sites of the habitat are similar.

 In this paper, we only consider the small advection and assume that
$$(H_2)\qquad \alpha<2\sqrt{d_I [\frac 1T\int^T_0(\beta_\infty(t)-\gamma_\infty(t))dt]}.$$
The initial distribution of the infected populations $I_0(x)$ is nonnegative and satisfies
\begin{equation}
I_0\in C^2[-h_0, h_0],\, I_0(-h_0)=I_0(h_0)=0\, \textrm{and} \ 0<I_0(x)\leq N^*,\, x\in (-h_0, h_0),
\label{Ae2}
\end{equation}
where the condition (\ref{Ae2}) indicates that at the beginning, the infected exists in the area with $x\in (-h_0, h_0)$,
but for the area $|x|\geq h_0$, no infected happens yet. Therefore,
the model means that beyond the left boundary $x=g(t)$ and the right boundary $x=h(t)$, there is only susceptible, no infected individuals.

The equations governing the free boundary, the spreading front,
$h'(t)=-\mu I_{x}(h(t),t)$ and $g'(t)=-\mu I_{x}(g(t),t)$, are the special cases of the well-known
Stefan condition, which has been established in \cite{LIN} for the diffusive populations.
The positive constant $\mu$ measures the ability of the infected transmit and diffuse towards the new area.

After we finished the first version of this paper, we found the papers (\cite{CLW, Wang15}) dealing with similar problem, which describing spatial spreading of the species without advection.
In our paper, we considered the effect of advection by using the basic reproduction numbers. We therefore emphasized the different consideration and omitted some similar proofs.

This rest of the paper is arranged as follows. In the next section, the global
existence and uniqueness of the solution to
(\ref{a3}) are presented by using a contraction mapping theorem,
comparison principle is also employed. Section 3 is devoted to introducing the basic reproduction
numbers and deriving their analytical properties. Sufficient conditions for the disease to vanish or spread are given in section 4.
The asymptotic spreading speeds are also presented.

\section{Preliminaries}

In this section, we first present some fundamental results on solutions of problem \eqref{a3}, we omit the proof since it is standard, see also Lemma 2.2,
Theorems 2.1 and 2.2 in \cite{GLZ}.

\begin{thm} For any given $I_0$ satisfying \eqref{Ae2}, and any $\nu \in (0, 1)$, problem \eqref{a3}
admits a unique global solution
$$(I; g, h)\in C^{ 1+\nu,(1+\nu)/2}([g(t), h(t)]\times [0,+\infty))\times C^{1+\nu/2}([0,T])\times C^{1+\nu/2}([0,+\infty));$$
moreover,
\[
0<I(x, t)\leq N^*\; \mbox{ for } g(t)<x<
h(t),\; t\in (0,T_0]. \]
\[
-C_1\leq g'(t)<0\; \mbox{and}\; 0<h'(t)\leq C_1 \;  \; t\in (0,T_0]. \]
for some constants $C_1$ and $T_0$.\label{exist}
\end{thm}

For later applications, we exhibit the comparison principle, which is similar to Lemma 3.5 in \cite{DL}.
\begin{lem} (The Comparison Principle)
  Assume that $T\in (0,\infty)$, $\overline g, \overline
h, \underline g, \underline
h$ $\in C^1([0,T])$, $\overline I(x,t)\in C([\overline g(t), \overline h(t)]\times [0, T])\cap
C^{2,1}((\overline g(t), \overline h(t))\times (0, T])$, $\underline I(x,t)\in C([\underline g(t), \underline h(t)]\times [0, T])\cap
C^{2,1}((\underline g(t), \underline h(t))\times (0, T])$, and
\begin{eqnarray*}
\left\{
\begin{array}{lll}
\overline I_{t}-d_I \overline I_{xx}+\alpha \overline I_x\geq (\beta (x,t) -\gamma(x,t))\overline I-\frac {\beta (x,t)}{N^*} \overline I^2,
&\overline g(t)<x<\overline h(t), \ 0<t\leq T,\\
\underline I_{t}-d_I \underline I_{xx}+\alpha \underline I_x\leq (\beta (x,t) -\gamma(x,t))\underline I-\frac {\beta (x,t)}{N^*} \underline I^2,
&\underline g(t)<x<\underline h(t), \ 0<t\leq T,\\
\overline I(\overline g(t), t)=0,\; \overline g'(t)\leq -\mu \overline I_x(\overline g(t), t),\quad & 0<t\leq T,\\
\underline I(\underline g(t), t)=0,\; \underline g'(t)\geq -\mu \underline I_x(\underline g(t), t),\quad & 0<t\leq T,\\
\overline I(\overline h(t), t)=0,\; \overline h'(t)\geq -\mu \overline I_x(\overline h(t), t),\quad & 0<t\leq T,\\
\underline I(\underline h(t), t)=0,\; \underline h'(t)\leq -\mu \underline I_x(\underline h(t), t),\quad & 0<t\leq T,\\
\overline g(0)\leq -h_0<h_0\leq \overline h(0), \, I_{0}(x)\leq \overline I(x, 0),\;  &-h_0\leq x\leq h_0,\\
-h_0\leq \underline g(0)\leq \underline h(0)\leq h_0, \, \underline I(x,0)\leq I_{0}(x),\;  &\underline g(0)\leq x\leq \underline h(0),
\end{array} \right.
\end{eqnarray*}
then the solution $(I(x,t); g(t), h(t))$ to the free boundary problem $(\ref{a3})$ satisfies
$$\overline g(t)\leq g(t)\leq\underline g(t),\ \underline h(t)\leq h(t)\leq\overline h(t)\, \mbox{for}\, t\in [0, T],$$
$$\underline I(x, t) \leq I(x, t)\, \mbox{for}\,  (x, t)\in [\underline g(t), \underline h(t)]\times [0, T],$$
$$I(x, t) \leq \overline I(x, t)\, \mbox{for}\,  (x, t)\in [g(t), h(t)]\times [0, T].$$
\end{lem}

\bigskip
The pair $(\overline u; \overline g, \overline h)$ in Lemma 2.2 is usually called an upper solution
of the problem \eqref{a3} and $(\underline u; \underline g, \underline h)$ is then called a lower solution.
To examine the dependence of the solution on the expanding capability $\mu$,
 we write the solution as $(I^{\mu}; g^{\mu}, h^{\mu})$. As a corollary of Lemma 2.2, we have the following monotonicity:

\begin{cor} For fixed $I_0, \alpha, h_0, \beta (x,t)$ and $\gamma (x,t)$.
If $\mu_1\leq \mu_2$. Then $I^{\mu_1}(x, t)\leq I^{\mu_2}(x, t)$ in $[g^{\mu_1}(t), h^{\mu_1}(t)]\times (0, \infty)$
 and $g^{\mu_2}(t)\leq g^{\mu_1}(t)$, $h^{\mu_1}(t)\leq h^{\mu_2}(t)$ in $(0, \infty)$.
 \label{mu}
\end{cor}

\section{Basic reproduction numbers}

In this section, we first present the basic reproduction number and its properties
for the corresponding system in a fixed interval, and then define the basic reproduction number
for the free boundary problem \eqref{a3}.
The basic reproduction numbers are related to the eigenvalues of corresponding periodic-parabolic eigenvalue problems.

Consider the reaction-diffusion-advection problem
\begin{eqnarray}
\left\{
\begin{array}{lll}
I_t-d_I I_{xx}+\alpha I_x=-\gamma(x,t)I,\; &x\in (h_1,h_2),\, t>0,  \\
\phi(h_1,t)=\phi(h_2,t)=0, & t>0,
\end{array} \right.
\label{v0}
\end{eqnarray}
Let $U(t,s)$ be the evolution operator of (\ref{v0}), then there exist constants $C, \omega>0$ such that
$||U(t,s)||\leq C{\textrm e}^{-\omega (t-s)}$ for any
$t\geq s, t, s\in \R$. As in \cite{PZ}, we introduce the next generation operator
$$L(\psi)(t):=\int^\infty_0 U(t,t-s)\beta(\cdot,t-s)\psi(\cdot,t-s)ds,$$
where  $\psi\in C_T:=\{u:\, t\rightarrow u(t)\in C[h_1, h_2],\, u(0)(x)=u(T)(x)\}.$ From the definition,
we know that $L$ is continuous, strongly positive and compact on $C_T$. We now define the basic reproduction number
of system (\ref{a3}) as the spectral radius of $L$, that is,
$$R^D_0=R^D_0(d_I, \alpha,  \beta(x,t), \gamma(x,t), [h_1, h_2])=\rho(L).$$

With the above definition, we have the following equivalent characterizations.
\begin{lem} $(i)$ $R_0^D=\mu_0$, where $\mu_0$ is the unique principal eigenvalue of periodic-parabolic eigenvalue problem
\begin{eqnarray}
\left\{
\begin{array}{lll}
\phi_t-d_I \phi_{xx}+\alpha \phi_x=-\gamma(x,t)\phi+\frac{\beta(x,t)}{\mu_0}\phi,\; &x\in (h_1,h_2),\, t>0,  \\
\phi(h_1,t)=\phi(h_2,t)=0, & t>0,\\
\phi(x,0)=\phi(x,T),&x\in [h_1,h_2].
\end{array} \right.
\label{v1}
\end{eqnarray}
$(ii)$ $1-R_0^D$ has the same sign as $\lambda_0$, where $\lambda_0$ is the principal eigenvalue of periodic-parabolic eigenvalue problem
\begin{eqnarray}
\left\{
\begin{array}{lll}
\psi_t-d_I \psi_{xx}+\alpha \psi_x=\beta(x,t) \psi-\gamma(x,t)\psi+\lambda_0 \psi,\; &
x\in (h_1,h_2),\, t>0,  \\
\psi(h_1,t)=\psi(h_2,t)=0,& t>0,\\
\psi(x,0)=\psi(x,T), &x\in [h_1, h_2].
\end{array} \right.
\label{B1f}
\end{eqnarray}
\end{lem}

The proof of this lemma is similar as that of Lemmas 2.1 and 2.2 in \cite{PZ}.
The existence of the unique principal eigenvalue $\mu_0$ of (\ref{v1}) and $\lambda_0$ of (\ref{B1f}) can be seen from Section 16 (Theorem 16.1) and Section 14  in \cite{Hess}, respectively. Moreover,  the eigenfunction $\phi(x,t)$ of
(\ref{v1}) corresponding to $\mu_0$ and  the eigenfunction $\psi(x,t)$ of (\ref{B1f}) corresponding to $\lambda_0$ are positive in $(h_1, h_2)\times \R$.

To see the properties of the basic reproduction number $R_0^D$, let us see the two special cases.
In the first case, if $\beta(x,t)$ and $\gamma(x,t)$ are spatially homogeneous, we have:
\begin{thm} If $\beta(x,t)\equiv \beta (t)$ and $\gamma (x,t)\equiv \gamma (t)$, then\\
$(i)$  $R_0^D(d_I, \alpha, \beta(t),\gamma(t), [h_1, h_2])=\frac {\frac 1T\int^T_0\beta(t)dt}{d_I(\frac {\pi}{h_2-h_1})^2+\frac{\alpha^2}{4d_I}+\frac 1T\int^T_0\gamma(t)dt}$;\\
$(ii)$  $R_0^D(d_I, \alpha, \beta(t),\gamma(t), [h_1, h_2])$ is monotone increasing with respect to $\beta(t)$, and decreasing with respect to $\alpha$ and $\gamma(t)$;\\
$(iii)$ $R_0^D(d_I,\alpha, \beta(t),\gamma(t),\cdot)$ is an increasing function for fixed $d_I, \alpha, \beta(t), \gamma(t)$, in the sense that $R_0^D(d_I,\alpha, \beta(t),\gamma(t),I_1)\leq R_0^D(d_I,\alpha, \beta(t),\gamma(t),I_2)$ provided that $I_1\subseteq I_2\subseteq \mathbb{R}^1$;\\
$(iv)$
$$\lim_{(h_2-h_1)\rightarrow0^+} R_0^D(d_I,\alpha, \beta(t), \gamma(t),[h_1,h_2])=0,$$ and $$\lim_{(h_2-h_1)\rightarrow+\infty}R_0^D(d_I,\alpha, \beta(t), \gamma(t),[h_1,h_2])=\frac{\frac{1}{T}\int_0^T\beta(t)dt}{\frac{\alpha^2}{4d_I}+\frac{1}{T}\int_0^T\gamma(t)dt}.$$ Moreover, there exists a positive constants $h^*$ such that $R_0^D(d_I,\alpha, \beta(t), \gamma(t),[h_1,h_1+h^*])=1$,  $R_0^D(d_I,\alpha, \beta(t), \gamma(t),[h_1, h_1+h])>1$ for $h>h^*$ and $R_0^D(d_I,\alpha, \beta(t), \gamma(t),[h_1, h_1+h])<1$ for $h<h^*$ provided that $(H_1)$ holds.
\end{thm}
\bpf Let
$$\mu^*=\frac {\frac 1T\int^T_0\beta(t)dt}{d_I(\frac {\pi}{h_2-h_1})^2+\frac{\alpha^2}{4d_I}+\frac 1T\int^T_0\gamma(t)dt},$$
$$u(x,t)=\textrm{e}^{\int^t_0(-\gamma(s)-d_I(\frac {\pi}{h_2-h_1})^2-\frac{\alpha^2}{4d_I}+\frac 1{\mu^*}\beta(s))ds}\textrm{e}^{\frac{\alpha x}{2d_I}}\cos \frac {\pi}{h_2-h_1}(x-\frac{h_2-h_1}2),$$
it is easy to see that $u(x,t)=u(x,t+T)$ and $u(x,t)$ is a positive $T-$periodic solution to (\ref{v1}) with $\mu_0=\mu^*$.
The result of $(i)$ is follows from the uniqueness of the principal eigenvalue of (\ref{v1}), and
the conclusions of $(ii)-(iv)$ follow from the expression of $R_0^D(d_I,\alpha, \beta(t), \gamma(t), [h_1, h_2])$ in $(i)$ directly.
\epf
\bigskip

For the other special case, $\beta(x,t)\equiv \beta (x)$ and $\gamma (x,t)\equiv \gamma (x)$, we have
\begin{thm} The following assertions hold.

$(a)$ $R_0^{D}(d_I, \alpha,  \beta(x), \gamma(x), [h_1, h_2])=\ \sup_{\psi\in H^1_0(h_1, h_2),\psi\neq 0} \{\frac{\int_{h_1}^{h_2} \beta
\psi^2dx}{\int_{h_1}^{h_2} (d_I\psi_x^2+\frac{\alpha^2}{4d_I}\psi^2+\gamma \psi^2)dx}\}$ is a positive and monotone decreasing function of $\alpha$;

$(b)$ If $I_1\subseteqq I_2 \subseteqq \mathbb{R}^1$, then $R_0^{D}(d_I,\alpha, \beta(x), \gamma(x), I_1)\leq R_0^{D}(d_I, \alpha, \beta(x), \gamma(x), I_2)$, with strict inequality if $I_2 \setminus I_1$ is a nonempty open set.

$(c)$ If $\beta (x,t)\equiv \beta_{\infty}$ and $\gamma(x,t)\equiv \gamma_{\infty}$, then
$$R_0^{D}(d_I, \alpha,  \beta_\infty, \gamma_\infty, [h_1, h_2])=\frac { \beta_{\infty}}{d_I(\frac {\pi}{h_2-h_1})^2+\frac{\alpha^2}{4d_I}+\gamma_{\infty} }.$$

$(d)$ If $(H_1)$ holds, then
$$\lim_{(h_2-h_1)\rightarrow0^+} R_0^D(d_I,\alpha, \beta(x), \gamma(x),[h_1,h_2])=0,$$ and $$\lim_{(h_2-h_1)\rightarrow+\infty}R_0^D(d_I,\alpha, \beta(x), \gamma(x), [h_1,h_2])\geq \frac{\beta_\infty}{\frac{\alpha^2}{4d_I}+\gamma_\infty}.$$
Furthermore, if the assumption $(H_2)$ holds, we can find a positive constant $h^*$ such that $R_0^{D}(d_I, \alpha,  \beta(x), \gamma(x), [h_1, h_1+h^*])=1$ and
$R_0^{D}(d_I, \alpha,  \beta(x), \gamma(x), [h_1, h_1+h])>1$ for $h>h^*$, $R_0^{D}(d_I, \alpha,  \beta(x), \gamma(x), [h_1, h_1+h])<1$ for $h<h^*.$
\label{basic2}
\end{thm}
\bpf
Let $\mu_0$ be the principal eigenvalue of the elliptic problem
\begin{eqnarray}
\left\{
\begin{array}{lll}
-d_I \phi_{xx}+\alpha \phi_x=-\gamma(x)\phi+\frac{\beta(x)}{\mu_0} \phi,\; &
x\in (h_1, h_2),  \\
\phi(h_1)=\phi(h_2)=0. &
\end{array} \right.
\label{B11f}
\end{eqnarray}
It follows from Lemma 3.1 (i) and the variational methods that
$$R_0^D(d_I, \alpha,  \beta(x), \gamma(x), [h_1, h_2])=\mu_0=\ \sup_{\phi\in H^1_0(h_1, h_2),\phi\neq 0} \{\frac{\int_{h_1}^{h_2} \beta
e^{-\alpha x/d_I}\phi^2dx}{\int_{h_1}^{h_2} (d_Ie^{-\alpha x/d_I}\phi_x^2+\gamma e^{-\alpha x/d_I}\phi^2)dx}\}.$$
If $\phi\in H^1_0(h_1, h_2)$, then $\psi=e^{-\alpha x/(2d_I)}\phi\in H^1_0(h_1, h_2)$,
therefore taking $\phi=e^{\alpha x/(2d_I)}\psi$ gives that
$$R_0^{D}(d_I, \alpha,  \beta(x), \gamma(x), [h_1, h_2])=\ \sup_{\psi\in H^1_0(h_1, h_2),\psi\neq 0} \{\frac{\int_{h_1}^{h_2} \beta
\psi^2dx}{\int_{h_1}^{h_2} (d_I\psi_x^2+\frac{\alpha^2}{4d_I}\psi^2+\gamma \psi^2)dx}\}.$$
Then (a)-(c) hold directly from the formulation of $R_0^D(d_I, \alpha,  \beta(x), \gamma(x), [h_1,h_2]$).

We then verify $(d)$. It's easy to see that
$$\lim_{(h_2-h_1)\rightarrow0^+} R_0^D(d_I,\alpha, \beta(t), \gamma(t),[h_1,h_2])=0.$$
Assertion $\lim_{(h_2-h_1)\rightarrow+\infty}R_0^D(d_I,\alpha, \beta(t), \gamma(t),[h_1,h_2])\geq \frac{\beta_\infty}{\frac{\alpha^2}{4d_I}+\gamma_\infty}$ follows from the proof of Theorem 3.2 (e) \cite{GLZ} by little modification. In view of (b) and $(H_2)$, we can find a constant $h^*>0$ such that $R_0^D([h_1, h_1+h^*]):=R_0^D(d_I,\alpha, \beta(x), \gamma(x), [h_1, h_1+h^*])=1$ and
$$R_0^D([h_1, h_1+h])>1\;\mbox{for}\; h>h^*, R_0^D([h_1, h_1+h])<1\;\mbox{for}\;h<h^*.$$
\epf

\bigskip

With the above properties in mind, we give some properties of $R^D_0(d_I, \alpha,  \beta(x,t), \gamma(x,t), [h_1, h_2])$. Firstly, we present some monotonicity.
\begin{thm} The following assertions hold.

$(1)$ $R^D_0(d_I, \alpha,  \beta(x,t), \gamma(x,t), [h_1, h_2])$ is a monotone increasing with respect to $\beta(x,t)$ and decreasing with respect to $\gamma(x,t)$, that is,
$$R^D_0(d_I, \alpha,  \beta_*(x,t), \gamma^*(x,t), [h_1, h_2])\leq R^D_0(d_I, \alpha,  \beta^*(x,t), \gamma_*(x,t), [h_1, h_2])$$
if $\beta_*(x,t)\leq \beta^*(x,t)$ and $\gamma_*(x,t)\leq \gamma^*(x,t)$;

$(2)$ $R^D_0(d_I, \alpha,  \beta(x,t), \gamma(x,t), [h_1, h_2])$ is a positive bounded function and satisfies
$$\frac{\frac 1T\int^T_0\beta_m(t)dt}{d_I(\frac {\pi}{h_2-h_1})^2+\frac{\alpha^2}{4d_I}+\frac 1T\int^T_0\gamma_M(t)dt} \leq
R_0^{D}\leq \frac {\frac 1T\int^T_0\beta_M(t)dt}{d_I(\frac {\pi}{h_2-h_1})^2+\frac{\alpha^2}{4d_I}+\frac 1T\int^T_0\gamma_m(t)dt},$$
where $\gamma_m(t)=\min_{x\in [h_1, h_2]} \gamma(x,t)$, $\gamma_M(t)=\max_{x\in [h_1, h_2]} \gamma(x,t)$
 and $\beta_m(t), \beta_M(t)$ are defined similarly;

$(3)$ $R_0^D(d_I,\alpha, \beta(x,t),\gamma(x,t),\cdot)$ is an increasing function for fixed $d_I, \alpha, \beta(x,t), \gamma(x,t)$, in the sense that $R_0^D(d_I,\alpha, \beta(x,t),\gamma(x,t),I_1)\leq R_0^D(d_I,\alpha, \beta(x,t),\gamma(x,t),I_2)$ provided that $I_1\subseteq I_2\subseteq \mathbb{R}^1$;

$(4)$ Assume that $(H_1)$ holds, then
$$\lim_{(h_2-h_1)\rightarrow0^+} R_0^D(d_I,\alpha, \beta(x,t), \gamma(x,t),[h_1,h_2])=0,$$  $$\liminf_{(h_2-h_1)\rightarrow+\infty}R_0^D(d_I,\alpha, \beta(x,t), \gamma(x,t),[h_1,h_2])\geq \frac{\frac{1}{T}\int_0^T\beta_\infty(t)dt}{\frac{\alpha^2}{4d_I}+\frac{1}{T}\int_0^T\gamma_\infty(t)dt}.$$
Furthermore, there exists a unique positive constant $h^*$ such that
$$R_0^D(h^*):=R_0^D(d_I,\alpha, \beta(x,t), \gamma(x,t),[h_1,h_1+h^*])=1$$
and $R_0^D([h_1, h_1+h)>1$  for $h>h^*$, $R_0^D([h_1, h_1+])<1$ for $h<h^*$ provided that $(H_2)$ holds.
\label{basic3}
\end{thm}
\bpf
$(1)$  For any $\beta_1(x,t)\leq\beta_2(x,t)$, denote $\mu_0^i (i=1,2)$ is the principal eigenvalue of (\ref
{v1}) and the corresponding eigenfunction is $\phi^i>0 (i=1,2)$ in $(h_1,h_2)$. It's well-known (Theorem 7.2  in \cite{Hess}) that $\mu_0^2$ is also the eigenvalue of
\begin{eqnarray}
\left\{
\begin{array}{lll}
-\psi_t-d_I \psi_{xx}-\alpha \psi_x=-\gamma(x,t)\phi+\frac{\beta_2(x,t)}{\mu_0^2} \psi,\; &x\in (h_1, h_2),\, t>0,  \\
\psi(h_1,t)=\psi(h_2,t)=0, & t>0,\\
\psi(x,0)=\psi(x,T),&x\in [h_1, h_2].
\end{array} \right.
\label{22}
\end{eqnarray}
It's corresponding eigenfunction is denoted by $\psi_2>0$ in $[h_1, h_2]$.

Next we apply the multiply-multiply-subtract-integrate trick.
Multiply the equation of $\phi_1$ by $\psi_2$ and the equation of $\psi_2$ by $\phi_1$. Then subtract the two equations and integrate over $(h_1, h_2)\times(0,T)$ obtain
$$ \frac{\int_0^T \int_{h_1}^{h_2} \beta_1(x,t)\phi_1\psi_2dxdt}{\mu_0^1}=\frac{\int_0^T\int_{h_1}^{h_2}\beta_2(x,t)\phi_1\psi_2dxdt}{\mu_0^2}.$$
Hence, $\mu_0^1\leq \mu_0^2$. From Lemma 3.1 (i), it follows that $$R^D_0(d_I,\alpha, \beta_1(x,t), \gamma(x,t),[h_1,h_2])\leq R^D_0(d_I,\alpha, \beta_2(x,t), \gamma(x,t),[h_1,h_2])$$
Therefore, $R^D_0(d_I,\alpha, \beta(x,t), \gamma(x,t),[h_1,h_2])$ is an increasing function of $\beta(x,t)$. Similarly, we deduce that $R^D_0(d_I,\alpha, \beta(x,t), \gamma(x,t),[h_1, h_2])$ is decreasing of $\gamma(x,t)$.

$(2)$  We first get an upper bound of $R^D_0$ by considering the eigenvalue problem
 \begin{eqnarray}
\left\{
\begin{array}{lll}
\tilde \phi_t-d_I \tilde \phi_{xx}+\alpha \tilde \phi=-\gamma_m(t)\tilde \phi+\frac{\beta_M(t)}{\tilde \mu} \tilde \phi,\; &x\in (h_1,h_2),\, t>0,  \\
\tilde \phi(h_1,t)=\tilde \phi(h_2,t)=0, & t>0,\\
\tilde \phi(x,0)=\tilde \phi(x,T),&x\in [h_1,h_2].
\end{array} \right.
\label{v13}
\end{eqnarray}
It follows from Theorem 3.2 that
$$\tilde \mu = \frac {\frac 1T\int^T_0\beta_M(t)dt}{d_I(\frac {\pi}{h_2-h_1})^2+\frac{\alpha^2}{4d_I}+\frac 1T\int^T_0\gamma_m(t)dt}.$$
We apply the monotonicity of $\beta(x,t)$, $\gamma(x,t)$ in (1) to deduce that
$$R^D_0(d_I,\alpha, \beta_1(x,t), \gamma(x,t),[h_1,h_2])\leq \tilde \mu.$$

Similarly, we obtain
$$R^D_0(d_I,\alpha, \beta_1(x,t), \gamma(x,t),[h_1,h_2])\geq \frac {\frac 1T\int^T_0\beta_m(t)dt}{d_I(\frac {\pi}{h_2-h_1})^2+\frac{\alpha^2}{4d_I}+\frac 1T\int^T_0\gamma_M(t)dt}.$$

$(3)$ Without loss of generality, we assume that $I_1=[0, l_1]$, $I_2=[0, l_2]$ and $0<l_1\leq l_2$.
We will show that $R_0^D(I_1) \leq R_0^D(I_2)$.
Denote $\phi$ be the eigenfunction corresponding to $R_0^D(I_1)$. Then we know that $R_0^D(I_2)$ is the principal eigenvalue of (\ref{22}) and the corresponding eigenfunction is $\psi(x,t)(>0)$.

Multiply the equation of $\psi$ by $\phi$ and the equation of $\phi$ by $\psi$,
then integrate them over $(0, l_1)\times(0, T)$ and subtract the resulting identities to obtain
$$\left( \frac{1}{R_0^D(I_1)}-\frac{1}{R_0^D(I_2)}\right)\int_0^T \int_0^{l_1}\beta(x,t)\phi \psi dxdt=-d_I\int_0^T\psi(l_1,t)\phi_x(l_1,t)dt.$$
Applying the Hopf Lemma to the equation of $\phi$ in $(0,l_1)\times(0,T)$ and we get $\phi_x(l_1,t)<0$. It follows that
$\left( \frac{1}{R_0^D(I_1)}-\frac{1}{R_0^D(I_2)}\right)\int_0^T \int_0^{l_1}\beta(x,t)\phi \psi dxdt>0$.
Therefore, $R_0^D(I_1) \leq R_0^D(I_2)$.

Finally, we verify $(4)$. Combine with the properties of $R_0^D(d_I, \alpha, \beta(x,t), \gamma(x,t), [h_1,h_2])$ in (2), we easily deduce that
$$\lim_{(h_2-h_1)\rightarrow +\infty}R_0^D(d_I, \alpha, \beta(x,t), \gamma(x,t), [h_1,h_2])=0.$$

Now, we show that
$$\liminf_{(h_2-h_1)\rightarrow+\infty}R_0^D(d_I,\alpha, \beta(x,t), \gamma(x,t),[h_1,h_2])\geq \frac {\frac{1}{T}\int_0^T\beta_{\infty}(t)dt}{\frac{\alpha^2}{4d_I}+\frac{1}{T}\int_0^T\gamma_{\infty}(t)dt}.$$
It follows from the assumption of $(H_1)$ that for any $\varepsilon>0$,
there exists a positive constant $L_0$ such that for $|x|\geq L_0$,
\begin{equation}\label{beta}
|\beta(x,t)-\beta_\infty|<\varepsilon,\ |\gamma(x,t)-\gamma_\infty|<\varepsilon.
\end{equation}
Without loss of generality, we assume that $h_1= 0$ and $h_2\to \infty$.
For any $L\geq 2L_0$, using the monotonicity of $R^D_0$ with respect to $\beta, \gamma$ and interval $[h_1, h_2]$ gives
\begin{eqnarray*}
& &R_0^D(d_I, \alpha,  \beta(x,t), \gamma(x,t), [0, L])\\
&\geq& R_0^D(d_I, \alpha,  \beta(x,t), \gamma(x,t), [L/2, L])\\
&\geq &R_0^D(d_I, \alpha,  \beta_\infty(t)-\varepsilon, \gamma_\infty(t)+\varepsilon, [L/2, L])\ (\textrm{ by }(\ref{beta}))\\
&=&\frac {\frac 1T\int^T_0(\beta_\infty(t)-\varepsilon)dt}{d_I(\frac {\pi}{L-L/2})^2+\frac{\alpha^2}{4d_I}+\frac 1T\int^T_0(\gamma_\infty(t)+\varepsilon)dt}\ (\textrm{ by Theorem 3.2 }),
\end{eqnarray*}
which together with the arbitrariness of small $\varepsilon$ gives
$$\liminf_{h_2\to +\infty}\, R_0^D([0,h_2])\geq \liminf_{L\to \infty}\, R_0^D([0,L])\geq\frac {\frac 1T\int^T_0\beta_\infty(t)dt}{\frac{\alpha^2}{4d_I}+\frac 1T\int^T_0\gamma_\infty(t)dt}.$$

In view of the assumption of $(H_2)$ holds, then we know that
 $$\liminf_{(h_2-h_1)\to +\infty}R_0^D(d_I, \alpha,  \beta(x,t), \gamma(x,t), [h_1, h_2])\geq 1.$$ As $R_0^D(d_I, \alpha,  \beta(x,t), \gamma(x,t), [h_1, h_2])$ is increasing with $(h_2-h_1)$, we can find $h^*>0$ such that $R_0^D([h_1, h_1+h^*]):=R_0^D(d_I, \alpha,  \beta(x,t), \gamma(x,t), [h_1, h_1+h^*])=1$ and
$$R_0^D([h_1, h_1+h])>1\;\mbox{for}\;h>h^*, R_0^D([h_1, h_1+h])<1\;\mbox{for}\;h<h^*.$$
\epf

Noticing that the interval $[g(\tau), h(\tau)]$, where the solution for the free boundary problem \eqref{a3} exist,
is changing with $\tau$, so the basic reproduction number is not a constant and
should be changing.
Now we introduced the basic reproduction number $R_0^F(\tau)$ for the free boundary problem \eqref{a3} by
$$R_0^F(\tau):=R^D_0(d_I, \alpha,  \beta(x,t), \gamma(x,t), [g(\tau), h(\tau)]),$$
it follows from Theorems \ref{exist} and \ref{basic3} that
\begin{thm} $R_0^F(\tau)$ is strictly monotone increasing function of $\tau$, that is if $\tau_1<\tau_2$, then $R_0^F(\tau_1)<R_0^F(\tau_2)$.
Moreover, if $(H_1)$ holds and $h(\tau)-g(\tau)\to \infty$ as $\tau\to \infty$, then $\lim_{\tau\to \infty} R_0^F(\tau)\geq \frac {\frac 1T\int^T_0(\beta_\infty(t)dt}{\frac{\alpha^2}{4d_I}+\frac 1T\int^T_0\gamma_\infty(t)dt}$.
\end{thm}

\begin{rmk} Assume that $(H_1)$ and $(H_2)$ hold. By Theorem 3.5, we have
 $R_0^F(\tau_0)>1$ for some $\tau_0>0$ provided that $h(\tau)-g(\tau)\to \infty$ as $\tau\to \infty$.\label{rem1}
\end{rmk}

\section{Spreading-vanishing}
It follows from Theorem \ref{exist} that $x=g(t)$ is monotonic decreasing and $x=h(t)$ is monotonic increasing, so
there exist $g_\infty\in [-\infty, -h_0)$ and $h_\infty\in (h_0, +\infty]$ such that $\lim_{t\to +\infty} \ g(t)=g_\infty$
and $\lim_{t\to +\infty} \ h(t)=h_\infty$. The following spreading-vanishing dichotomy has been given in \cite{CLW, Wang15} for the free boundary problem in time-periodic environment without advection.
\begin{lem}  Assume that $(H_1)$ and $(H_2)$ hold. Then, the following alternative
holds:\\
Either $(i)$ vanishing: $-\infty<g_\infty<h_\infty<\infty$, and
$$
R_0^{D}(d_I, \alpha,  \beta(x,t), \gamma(x,t), [g_\infty, h_\infty])\leq 1\ \textrm{and}\
\lim_{t\to\infty} \|I(\cdot, t)\|_{C([g(t),\,h(t)])}=0;$$
or $(ii)$ spreading: $-g_\infty=\infty=h_\infty$, then
$$\lim_{n\to \infty}I(x,t+nT)=\hat U(x,t)\, \textrm{locally uniformly for}\, (x,t)\in [0,\infty)\times [0,T],$$
where $\hat U(x,t)$ is the unique positive T-periodic solution of the problem
\begin{eqnarray}
\left\{
\begin{array}{lll}
U_t-d_I U_{xx}+\alpha U_x=(\beta (x,t) -\gamma(x,t))U-\frac {\beta (x,t)}{N^*} U^2,\; &x\in \R,\, 0\leq t\leq T,  \\
U(x,0)=U(x,T),&x\in \R.
\end{array} \right.
\label{v54}
\end{eqnarray}
\label{spr}
\end{lem}
\bpf  We first show that both $h_\infty$ and $g_\infty$ are finite or infinite simultaneously.
 In fact, if $h_\infty<\infty$, we can prove that $R^D_0(d_I, \alpha,  \beta(x,t), \gamma(x,t), [g_\infty, h_\infty])\leq 1$ as in Lemma 4.1 in \cite{GLZ} by contradiction, which together with the assumption, implies that $g_\infty>-\infty$ by Theorem 3.5.

In the case $-\infty<g_\infty<h_\infty<\infty$, since that $R^D_0(d_I, \alpha,  \beta(x,t), \gamma(x,t), [g_\infty, h_\infty])\leq 1$, we show that $\lim_{t\to
+\infty} \ ||I(\cdot,t)||_{C([g(t),\, h(t)])}=0$.

In fact, let $\overline I(x,t)$ denote the unique solution of the problem
\begin{eqnarray}
\left\{
\begin{array}{lll}
\overline I_t -d_I \overline I_{xx}+\alpha \overline I_x=\overline I(\beta(x,t)-\gamma(x,t))-\frac {\beta(x,t)}{N^*}\overline I^2,\;
&g_\infty<x<h_\infty,\, t>0, \\
\overline I(g_\infty, t)=0, \quad \overline I(h_\infty, t)=0,  & t>0,\\
\overline I(x,0)=\tilde I_0(x), &g_\infty<x<h_\infty,
\end{array} \right.
\label{m23}
\end{eqnarray}
with
\begin{eqnarray*}
\tilde I_0(x)= \left\{
\begin{array}{lll}
I_0(x),\ &g_0\leq x\leq h_0, \\
0, \ & \mbox{ otherwise}.
\end{array} \right.
\end{eqnarray*}
The comparison principle gives $0\leq I(x,t)\leq \overline I(x,t)$
for  $x\in [g(t), h(t)]$ and $t\geq 0$.

Using the fact $R^D_0(d_I, \alpha,  \beta(x,t), \gamma(x,t), [g_\infty, h_\infty])\leq 1$, we know that (Lemma 3.1 (ii))
the principal eigenvalue $\lambda_0$ of (3.3) is nonnegative.
the corresponding T-periodic problem
\begin{eqnarray}
\left\{
\begin{array}{lll}
U_t-d_I U_{xx}+\alpha U_x=U(\beta(x,t)-\gamma(x,t))-\frac {\beta(x,t)}{N^*}U^2,
\; &x\in (g_\infty, h_\infty),\, t>0,  \\
U(g_\infty,t)= U(h_\infty, t)=0, & t>0,\\
U(x,0)=U(x,T),&x\in [g_\infty, h_\infty].
\end{array} \right.
\label{v111}
\end{eqnarray}
admits only trivial solution $0$.
It is shown in \cite{Hess, Pao}, by the method of upper and lower solutions and
its associated monotone iterations, that the time-dependent solution $\overline I(x,t)$ converges to $0$ uniformly
in $[g_\infty,h_\infty]$ as $t\to \infty$.
Therefore, $\lim_{t\to +\infty} \ ||I(\cdot,t)||_{C([g(t),\, h(t)])}=0$.

In that case $-g_\infty=\infty=h_\infty$, we first consider problem \eqref{v54}.
When $\alpha=0$, the existence and uniqueness of positive T-periodic solution $U(x,t)$ of problem \eqref{v54}
  is directed from Theorem 1.3 in \cite{PD}.
 When $\alpha\neq 0$, the result still holds since $\tilde U:=N^*$ and $\hat U:=\varepsilon \phi_\delta(x)$ are the ordered upper and lower solutions of problem \eqref{v54}, where $\varepsilon$ and $\delta$ is a sufficiently small positive constant and $\phi_\delta$ satisfies
  \begin{eqnarray}
\left\{
\begin{array}{lll}
\phi_t-d_I \phi_{xx}+\alpha \phi_x=(\beta(x,t)-\gamma(x,t)-\delta)\phi,\; &x\in \R, \, t>0,  \\
\phi(x,0)=\phi(x,T),&x\in \R.
\end{array} \right.
\label{va1}
\end{eqnarray}
By the assumption that $(H_1)$ and $(H_2)$ hold, problem \eqref{va1} admits at one positive solution if $\delta<\frac 1T\int^T_0(\beta_\infty(t)-\gamma_\infty(t))dt-\frac{\alpha^2}{4d_I}$.

As to the limit, that is,
$$\lim_{n\to \infty}I(x,t+nT)=\hat U(x,t)\, \textrm{locally uniformly for}\, (x,t)\in [0,\infty)\times [0,T],$$
the proof is based on the upper and lower solutions method, see Lemma 4.2 in \cite{CLW} or Theorem 4.3 in \cite{Wang15}.
\epf

Now we give sufficient conditions so that the disease is spreading.
\begin{thm} If $R_0^F(t_0)\geq 1$ for some $t_0\geq 0$, then spreading must happen.
\end{thm}
\bpf Owing to $R_0^F(t_0)\geq 1$ for $t_0\geq 0$, $R_0^F(t_1)\geq 1$ for any $t_1>t_0$ by
the monotonicity in Theorem 3.4.

 In this case, we have that the periodic-parabolic eigenvalue problem
\begin{eqnarray}
\left\{
\begin{array}{lll}
\psi_t-d_I \psi_{xx}+\alpha \phi_x=\beta(x,t) \psi-\gamma(x,t)\psi+\lambda_0 \psi,\; &
x\in (g(t_1), h(t_1)),  \\
\psi(g(t_1),t)=\psi(h(t_1),t)=0,& t>0,\\
\psi(x,0)=\psi(x,T), &x\in [g(t_1), h(t_1)].
\end{array} \right.
\label{B1f1}
\end{eqnarray}
admits a positive solution $\psi(x)$ with $||\psi||_{L^\infty}=1$, where $\lambda_0$ is the principal eigenvalue. It follows from Lemma 3.1 that $\lambda_0<0$.

We are going to construct a suitable lower solutions to
\eqref{a3} and we define
$$
\underline{I}(x,t)=
\delta \psi(x,t),\quad g(t_1)\leq x\leq h(t_1), \, t\geq t_1,
$$
where $\delta$ is sufficiently small such that $0<\delta\leq \frac { N^*(-\lambda_0)}{\overline \beta}$,
$\delta \psi (x,t_1)\leq I(x, t_1)$ in $[g(t_1), h(t_1)]$ and $\overline \beta=\max_{[g(t_1),h(t_1)]\times[0,T]}\beta(x,t)$.

 Direct computations yield
\begin{eqnarray*}
& &\underline{I}_t-d_{I} \underline{I}_{xx}+\alpha \underline{I}_{x}-(\beta(x,t)-\gamma(x,t))\underline{I}+\frac{\beta(x,t)}{N^*}\underline I^{2}\\
& &=\delta \psi(x)[\lambda_0+\frac{\beta(x,t)}{N^*}\delta \psi(x)]\\
& &\leq 0
\end{eqnarray*}
for all $t>0$ and $g(t_1)<x<h(t_1)$.
Then we have
\begin{eqnarray*}
\left\{
\begin{array}{lll}
\underline{I}_t-d_{I}\underline{I}_{xx}+\alpha \underline{I}_x\leq (\beta(x,t)-\gamma(x,t))\underline{I}-\frac{\beta(x,t)}{N^*}\underline I^{2},\; &g(t_1)<x<h(t_1),\ t>t_1, \\
\underline{I}(g(t_1),t)=\underline{I}(h(t_1),t)=0,\; &t>t_1, \\
\underline{I}(x,t_1)\leq I(x, t_1)(x),\; &g(t_1)\leq x\leq h(t_1).
\end{array} \right.
\end{eqnarray*}
Using the comparison principle in the fixed interval $[g(t_1), h(t_1)]$ yields that  $I(x,t)\geq\underline I(x,t)$
in $[g(t_1), h(t_1)]\times [t_1,\infty)$. It follows that $\liminf_{t\to
+\infty} \ ||I(\cdot, t)||_{C([g(t_1), h(t_1)])}\geq \delta \psi(0)>0$ and therefore $h_\infty-g_\infty=+\infty$ by Lemma 4.1.
\epf

\bigskip
Similarly, we can also construct a suitable lower solution for $I$ to obtain sufficient conditions so that the disease is spreading and construct some suitable upper solutions to derive sufficient conditions so that the disease is vanishing. The proof of the next theorem will be omitted since it is an analogue of Lemma 3.7 in \cite{DL} or Lemma 2.8 in \cite{DG}.
\begin{thm} Suppose $R_0^F(0)(:=R_0^{D}(d_I, \alpha, \beta(x,t), \gamma (x,t), [-h_0, h_0])<1$.  Then
spreading happens
 if $\mu$ is sufficiently large; and vanishing happens
 if $||I_0(x)||_{C([-h_0, h_0])}$ or $\mu$ is sufficiently small.
 \label{condi}
\end{thm}

The following result is a consequence of Corollary \ref{mu} and Theorem \ref{condi}.
\begin{thm} (Sharp threshold)
Fixed $g_0, h_0$ and $I_0$. There exists $\mu^* \in [0, \infty)$
 such that spreading happens when $\mu> \mu^*$, and vanishing happens when $0<\mu\leq \mu^*$.
\end{thm}

Finally we give the asymptotic spreading speeds when the spreading happens.
\begin{thm}\label{spreadsp1} Assume that $(H_1)$ and $(H_2)$ hold. If $h_\infty=-g_\infty=+\infty$, then
\begin{align}
\lim_{t\to +\infty} \frac {h(t)}t=\frac 1T\int^T_0 k^*(\alpha, a, b)(t)dt,\quad \lim_{t\to +\infty} \frac {-g(t)}t=\frac 1T\int^T_0 k^*(-\alpha, a, b)(t)dt,
\label{limt}
\end{align}
\label{m10}
where $a(t)=\beta_\infty(t)-\gamma_\infty(t)$ and $b(t)=\frac{\beta_\infty(t)}{N^*}$, $(k,q) = (k^*(\alpha,a, b)(t), q^{*}(x,t))$ is the unique positive $T-$periodic solution of the problem
\begin{equation}\label{prob-q1}
\left\{
\begin{array}{ll}
q_t-d_I q'' +(k-\alpha)q'= q[a(t)-b(t) q],\, &x\in (0,\infty),\, t\in [0,T],\\
q(0,t)=0, \ & t\in [0,T],\\
q(x,0)=q(x,T),\ &x\in (0,\infty),\\
\mu q'(0,t)=k(t),\ &t\in [0,T],
\end{array}
\right.
\end{equation}
Moreover,
\begin{align}
0 < k^*(-\alpha, a, b) <k^*(0, a, b) <k^*(\alpha, a, b),
\label{daxiao}
\end{align}
 which means that the left boundary moves slower than the normal one without advection and the right boundary moves faster.
\end{thm}
\bpf
The existence and uniqueness of the solution $(k^*, q^{*})$ to problem \eqref{prob-q1} with $\alpha=0$ is given by Theorem 2.4 in \cite{DGP}. Checking the proof in \cite{DGP}, we found it still hold for $\alpha<2\sqrt{d_I(\overline {\beta_\infty}
-\overline {\gamma_\infty})}$. Moreover,
$$0<\overline {k^*}(\alpha, a, b)<2\sqrt{d_I(\overline {\beta_\infty}
-\overline {\gamma_\infty})}+\alpha,\ 0<\overline {k^*}(-\alpha, a, b)<2\sqrt{d_I(\overline {\beta_\infty}
-\overline {\gamma_\infty})}-\alpha,$$
where $\overline {\beta_\infty}=\frac 1 T\int^T_0 \beta_\infty(t)dt$ and $\overline {k^*}(\alpha, a, b), \overline {\gamma_\infty}$ are defined similarly.

The proof of \eqref{limt} is similar as that of Theorem 4.4 in \cite{DGP} with obvious modification,
 see also Theorem 5.5 in \cite{Wang15} or Corollary 3 in \cite{CLW}.

 \eqref{daxiao} can be established as Proposition 1.2 in \cite{GLL}.
\epf

\end{document}